\documentclass{article}
\usepackage{enumitem} 
\usepackage{arxiv}
\usepackage{subcaption}
\usepackage{amsmath}
\usepackage{amsthm}
\usepackage{amssymb}
\usepackage[T1]{fontenc}    
\usepackage{hyperref}       
\usepackage{url}            
\usepackage{booktabs}       
\usepackage{amsfonts}       
\usepackage{nicefrac}       
\usepackage{microtype}      
\usepackage{cleveref}       
\usepackage{lipsum}         
\usepackage{graphicx}
\usepackage{doi}
\usepackage[toc,page]{appendix}
\usepackage{xcolor}
\usepackage{etoolbox}
\usepackage{tikz}
\usepackage{esint}

\newtheoremstyle{mytheoremstyle} 
  {10pt}                        
  {10pt}                        
  {\normalfont}                  
  {}                             
  {\bfseries}                    
  {.}                            
  {5pt plus 1pt minus 1pt}        
  {}                             

\theoremstyle{mytheoremstyle}      

\setlist[enumerate,1]{label=\textup{\roman*)}}

\title{A Simple Proof of the Riemann Hypothesis}

\author{ Hatem A.~Fayed \\
	Applied Mathematics Program, Zewail City of Science, Technology and Innovation\\
October Gardens, 6th of
October, Giza 12578, Egypt\\
	\texttt{hfayed@zewailcity.edu.eg} \\
}

\renewcommand{\headeright}{}
\renewcommand{\undertitle}{}

\hypersetup{
pdftitle={Proof of the Riemann Hypothesis},
pdfsubject={Riemann Hypothesis},
pdfauthor={Hatem A.~Fayed},
pdfkeywords={Riemann zeta function, Riemann hypothesis, Non-trivial zeros, Critical line},
}

\begin{document}
\maketitle

\begin{abstract}
In this article, it is proved that the non-trivial zeros of the Riemann zeta function must lie on the critical line, known as the Riemann hypothesis.
\end{abstract}

\section{Riemann Zeta function}
The Riemann zeta function is defined over the complex plane as \cite{Riemann1859},
\begin{equation}
\label{eq1}
\zeta(s)=\sum_{n=1}^{\infty} \frac{1}{n^{s}} \quad, \Re(s)>1
\end{equation}
where $\Re(s)$ denotes the real part of $s$. There are several forms can be used for an analytic continuation for $\Re(s)>0$ such as \cite{Riemann1859,olver2010},
\begin{equation}
\label{eq31}
\zeta(s)=\sum_{n=1}^N \frac{1}{n^s}-\frac{N^{1-s}}{1-s}-s \int_N^{\infty} \frac{x-\lfloor x\rfloor}{x^{s+1}} d x, N\in\mathbb{N}
\end{equation}
where $\lfloor x\rfloor$ is the floor or integer part such that $x-1 < \lfloor x\rfloor \le x$ for real $x$.

and
\begin{equation}
\label{eq2}
\zeta(s)= \frac{1}{1-2^{1-s}}\eta(s),\qquad \eta(s)=\sum_{n=1}^{\infty} \frac{(-1)^{n+1}}{n^{s}}
\end{equation}
where $s\neq1+\frac{2 \pi k i}{\log (2)}, k=0, \pm 1, \pm 2 \ldots$ and $\eta(s)$ is the Dirichlet eta function (sometimes called the alternating zeta function).

Using Euler-Maclaurin formula, it can also be written as \cite{Riemann1859,olver2010},
\begin{equation}
\label{eq4}
\zeta(s)=\sum_{n=1}^N \frac{1}{n^s}-\frac{N^{1-s}}{1-s}-\frac{1}{2 N^s}+\sum_{i=1}^m\binom{s+2 i-2}{2 i-1} \frac{B_{2 i} }{2 i}N^{1-s-2 i}-\epsilon(m,N,s), \quad \Re(s)>-2 m ; m, N \in \mathbb{N}
\end{equation}
where
\begin{equation}
\begin{aligned}
\label{eq5}
\epsilon(m,N,s)&=\binom{s+2 m}{2 m+1}\int_{N}^{\infty}\frac{\bar{B}_{2m+1}(x)}{x^{s+2m+1}}dx,\\
\end{aligned}
\end{equation}
$B_{k}$ is the $k$-th Bernoulli number defined implicitly by,
\begin{equation}
\begin{aligned}
\label{eq6}
\frac{t}{e^t-1}=\sum_{k=0}^\infty B_k\frac{t^k}{k!},
\end{aligned}
\end{equation}
$B_k(x)$ is the $k$-th Bernoulli polynomial defined as the unique polynomial
of degree $k$ with the property that,
\begin{equation}
\begin{aligned}
\label{eq7}
\int_t^{t+1}B_k(x)dx=t^k,
\end{aligned}
\end{equation}
$\bar{B}_k(x)$ is the periodic function $B_k(x-\lfloor x\rfloor)$.

Riemann also extended $\zeta(s)$ to $\mathbb{C}$ as a meromorphic function, with only a simple pole at $s =1$ with residue $1$, by the functional equation,
\begin{equation}
\begin{aligned}
\label{eq72}
\zeta({s})=2^s\pi^{s-1}\sin\left(\frac{\pi s}{2}\right)\Gamma\left({1-s}\right)\zeta({1-s}), \qquad s\in \mathbb{C}\setminus \{0,1\}
\end{aligned}
\end{equation}
\section{Zeros of the Riemann Zeta Function}
The trivial zeros of the Riemann zeta function occur at the negative even integers; that is, $\zeta(-2n)=0,n\in\mathbb{N}$ \cite{Riemann1859}. On the other hand, the non-trivial zeros lie in the critical strip, $0\le \Re(s)\le1$. Both Hadamard \cite{Hadamard1896} and de la Vallee Poussin \cite{Vallee-Poussin1896} independently proved that there are no zeros on the boundaries of the critical strip (i.e. $\Re(s) = 0$ or $\Re(s) = 1$).
Gourdon and Demichel \cite{Gourdon2004} verified the Riemann Hypothesis until the $10^{13}$-th zero.

Mossinghoff and Trudgian \cite{Mossinghoff2015} proved that there are no zeros for $\zeta(\sigma+it)$ for $|t|\ge2$ in the region,
\begin{equation}
\label{eq8}
\sigma\ge1-\frac{1}{5.573412\log|t|}
\end{equation}
This represents the largest known zero-free region for the zeta-function within the critical strip for $3.06\times10^{10}<|t|<\exp(10151.5)\approx5.5\times10^{4408}$.

Due to the functional equation (\ref{eq72}) and the complex conjugation properties, the non-trivial zeros of $\zeta(s)$ are symmetric with respect to both the critical line and the real axis, that is $\zeta(s)=\zeta(1-s)=\zeta(\bar{s})=\zeta(1-\bar{s})=0$. According to equation (\ref{eq2}), the zeros of the Dirichlet eta function include all the non-trivial zeros of the zeta function. So, if $s$ is a non-trivial zero of the zeta function, then $\zeta(s)=\zeta(1-s)=\zeta(\bar{s})=\zeta(1-\bar{s})=0$. 

\section{Riemann Hypothesis}
All the non-trivial zeros of the Riemann zeta function lie on the critical line $\Re(s)=1/2$.
\begin{proof}
Assume that $s=\sigma+it, 0<\sigma<1, t\in \mathbb{R}, N\in\mathbb{N}, u\in \mathbb{R^+}$ such that $u<N$. In addition, let us restrict $u$ such that belongs to the truncated \ domain,
\begin{equation}
\begin{aligned}
\mathcal D=\bigcup_{j=j_0}^{\infty}D_j
, \qquad D_j=(a_j,b_j), \quad a_j=j+h/j, b_j=j+1-h/j,\quad 0<\delta<h<1, (2+|t|)^4=o(j_0)
\end{aligned}
\end{equation}
Note that $j=\lfloor u\rfloor$.

For $a,b\in \mathbb{R}$ such that $a<b$, the general Euler-Maclaurin formula \cite{Montgomery2006},\cite{ivic2003} is given by,
\begin{equation}
\begin{aligned}
\label{eq101}
\sum_{a < n \le b} f(n) = \int_a^b f(x) \, dx + \sum_{r=1}^{k} \frac{(-1)^r}{r!} \left[ \bar{B}_r(b) f^{(r-1)}(b) - \bar{B}_r(a) f^{(r-1)}(a) \right] - \frac{(-1)^{k}}{k!} \int_a^b \bar{B}_{k}(x) f^{(k)}(x) \, dx
\end{aligned}
\end{equation}
Substituting by $f(x)=1/x^s, a=u, b=N, k=4$,
\begin{equation}
\begin{aligned}
\label{eq102}
\sum_{u<n\le N}\frac{1}{n^s}&=\frac{N^{1-s}-u^{1-s}}{1-s}-\left[\frac{\bar{B}_1(N)}{N^{s}}-\frac{\bar{B}_1(u)}{u^{s}}\right]-\frac{s}{2!}\left[\frac{\bar{B}_2(N)}{N^{s+1}}-\frac{\bar{B}_2(u)}{u^{s+1}}\right]-\frac{s(s+1)}{3!}\left[\frac{\bar{B}_3(N)}{N^{s+2}}-\frac{\bar{B}_3(u)}{u^{s+2}}\right]\\
&-\frac{s(s+1)(s+2)}{4!}\left[\frac{\bar{B}_4(N)}{N^{s+3}}-\frac{\bar{B}_4(u)}{u^{s+3}}\right]+E(u,N,s)
\end{aligned}
\end{equation}
where 
\begin{equation}
\begin{aligned}
\label{eq103}
E(u,N,s)&=-\frac{s(s+1)(s+2)(s+3)}{4!}\int_{u}^{N} \frac{\bar{B}_4(x)}{x^{s+4}}\,dx,\\
B_1(x)&=x-\frac{1}{2},\quad B_2(x)=x^2-x+\frac{1}{6},\quad\frac{dB_n(x)}{dx}=nB_{n-1}(x),~n\in\mathbb{N}
\end{aligned}
\end{equation}
By adding $\sum_{1\le n\le u}\frac{1}{n^s}$ to both sides of equation (\ref{eq102}), we get
\begin{equation}
\begin{aligned}
\label{eq104}
\sum_{n=1}^{N}\frac{1}{n^s}&=\sum_{n\le u}\frac{1}{n^s}+\frac{N^{1-s}-u^{1-s}}{1-s}-\left[\frac{\bar{B}_1(N)}{N^{s}}-\frac{\bar{B}_1(u)}{u^{s}}\right]-\frac{s}{2!}\left[\frac{\bar{B}_2(N)}{N^{s+1}}-\frac{\bar{B}_2(u)}{u^{s+1}}\right]-\frac{s(s+1)}{3!}\left[\frac{\bar{B}_3(N)}{N^{s+2}}-\frac{\bar{B}_3(u)}{u^{s+2}}\right]\\
&-\frac{s(s+1)(s+2)}{4!}\left[\frac{\bar{B}_4(N)}{N^{s+3}}-\frac{\bar{B}_4(u)}{u^{s+3}}\right]+E(u,N,s)
\end{aligned}
\end{equation}
As $N\to\infty$, from equation (\ref{eq31}), we have,
\begin{equation}
\begin{aligned}
\label{eq105}
\sum_{n=1}^N\frac{1}{n^s} - \frac{N^{1-s}}{1-s}\to \zeta(s)
\end{aligned}
\end{equation}
and since $\left|\bar{B}_4\right|\le \bar{B}_4^{max}$, where $\bar{B}_4^{max}=1/30$, then the remainder in equation (\ref{eq104}) can be bounded by, 
\begin{equation}
\begin{aligned}
\label{eq106}
|E(u,\infty,s)|&\le \frac{\left|s(s+1)(s+2)(s+3)\right|\bar{B}_4^{max}}{4!}\int_{u}^{\infty} \frac{1}{x^{\sigma+4}}\,dx&= \frac{\left|s(s+1)(s+2)(s+3)\right|\bar{B}_4^{max}}{4!(\sigma+3)u^{\sigma+3}}
\end{aligned}
\end{equation}
Substituting by equations (\ref{eq105}) and (\ref{eq106}) in equation (\ref{eq104}) as $N\to\infty$, we get
\begin{equation}
\begin{aligned}
\label{eq107}
\sum_{n=1}^{j}\frac{1}{n^s}&=\zeta(s)+\frac{u^{1-s}}{1-s}-\frac{\bar{B}_1(u)}{u^{s}}-\frac{s\bar{B}_2(u)}{2u^{s+1}}-\frac{s(s+1)\bar{B}_3(u)}{3!u^{s+2}}+O\left(\frac{\left(1+|t|\right)^4}{u^{\sigma+3}}\right)
\end{aligned}
\end{equation}
An integral representation of the above sum can be obtained using the truncated Perron's formula \cite{Apostol1976} as follows.

For $z=c+iy$, $c=2$, $y\in\mathbb{R}$,
\begin{equation}
\begin{aligned}
\label{eq109}
\sum_{n=1}^\infty\frac{1}{n^{s+z}}=\sum_{n=1}^\infty\frac{1}{n^{2+\sigma+i(t+y)}}
\end{aligned}
\end{equation}
is absolutely and uniformly convergent.

Let $T=u^8$ and consider the following integral,
\begin{equation}
\begin{aligned}
\label{eq110}
\frac{1}{2\pi i} \int_{c-iT}^{c+iT}\frac{u^z}{z}\zeta(s+z)dz&=\frac{1}{2\pi i} \int_{c-iT}^{c+iT}\frac{u^z}{z}\sum_{n=1}^\infty\frac{1}{n^{s+z}}\,dz \\
&=\frac{1}{2\pi i} \int_{c-iT}^{c+iT}\frac{u^z}{z}\sum_{n=1}^{j}\frac{1}{n^{s+z}}\, dz +\frac{1}{2\pi i} \int_{c-iT}^{c+iT}\frac{u^z}{z}\sum_{n={j}+1}^\infty\frac{1}{n^{s+z}}\, dz \\
&=\frac{1}{2\pi i} \sum_{n=1}^{j}\frac{1}{n^s}\int_{c-iT}^{c+iT}\left(\frac{u}{n}\right)^z \frac{dz}{z} +\frac{1}{2\pi i} \sum_{n=j+1}^{\infty}\frac{1}{n^s}\int_{c-iT}^{c+iT}\left(\frac{u}{n}\right)^z \frac{dz}{z} 
\end{aligned}
\end{equation}
According to \cite{Apostol1976},
\begin{equation}
\label{eq111}
\begin{aligned}
\left|\frac{1}{2\pi i} \int_{c-iT}^{c+iT}\frac{a^z}{z}\,dz\right| & \le\frac{a^c}{\pi T \log(1/a)} && \text { if } 0<a<1\\
\left|\frac{1}{2\pi i} \int_{c-iT}^{c+iT}\frac{a^z}{z}\,dz-1\right| & \le\frac{a^c}{\pi T \log(a)} && \text { if } a>1\\
\end{aligned}
\end{equation}
So, for $n\ge{j}+1$, $n>u$, $a=u/n<1$ and $n/u>(u+h/j)/u$, then we have
\begin{equation}
\begin{aligned}
\log\left(\frac{u+h/j}{u}\right)
= \log\left(1+\frac{h}{ju}\right)
\sim \frac{h}{ju} \text{ for large } u
\end{aligned}
\end{equation}
Hence
\begin{equation}
\begin{aligned}
\label{eq112}
\left|\frac{1}{2\pi i} \sum_{n=j+1}^{\infty}\frac{1}{n^s}\int_{c-iT}^{c+iT}\left(\frac{u}{n}\right)^z \frac{dz}{z} \right|&\le \sum_{n={j}+1}^\infty\frac{1}{n^{\sigma}}\left[\frac{(u/n)^c}{\pi T \log (n/u)}\right]< \sum_{n={j}+1}^\infty\frac{1}{n^{\sigma}}\left[\frac{(u/n)^c}{\pi T \log \left(\frac{u+h/j}{u}\right)}\right]\\
&\lesssim\frac{ju^{c+1}}{\pi T h}\sum_{n={j}+1}^\infty\frac{1}{n^{\sigma+c}}= O\left(\frac{u^{c+2}}{T}\right)=O\left(\frac{1}{u^{4}}\right)
\end{aligned}
\end{equation}
On the other hand,
\begin{equation}
\begin{aligned}
\label{eq113}
\frac{1}{2\pi i} \sum_{n=1}^{j}\frac{1}{n^s}\int_{c-iT}^{c+iT}\left(\frac{u}{n}\right)^z \frac{dz}{z}-
\sum_{n=1}^{j}\frac{1}{n^s}&=\sum_{n=1}^{j}\frac{1}{n^s}\left[\frac{1}{2\pi i}\int_{c-iT}^{c+iT}\left(\frac{u}{n}\right)^z \frac{dz}{z}-1\right]
\end{aligned}
\end{equation}
For $1\le n\le{j}$, $a=u/n>1$ and $u/n\ge u/j>(j+h/j)/j$, then we have
\begin{equation}
\begin{aligned}
\label{eq114}
\left|\sum_{n=1}^{j}\frac{1}{n^s}\left[\frac{1}{2\pi i}\int_{c-iT}^{c+iT}\left(\frac{u}{n}\right)^z \frac{dz}{z}-1\right]\right|
&< \sum_{n=1}^{j}\frac{1}{n^{\sigma}}\left[\frac{(u/n)^c}{\pi T \log \left(\frac{j+h/j}{j}\right)}\right]\lesssim\frac{j^2 u^{c}}{\pi  T h}\sum_{n=1}^{j}\frac{1}{n^{\sigma+c}}\\
&\lesssim\frac{u^{c+2}}{\pi T h}\sum_{n=1}^{j}\frac{1}{n^{\sigma+c}}= O\left(\frac{1}{u^{4}}\right)
\end{aligned}
\end{equation}
From equations (\ref{eq110}), (\ref{eq112}), (\ref{eq113}) and (\ref{eq114}), we have
\begin{equation}
\begin{aligned}
\label{eq115}
\frac{1}{2\pi i} \int_{2-iT}^{2+iT}\frac{u^z}{z}\zeta(s+z)dz&=\sum_{n=1}^{j}\frac{1}{n^s}+O\left(\frac{1}{u^{4}}\right)
\end{aligned}
\end{equation}
Let
\begin{equation}
\begin{aligned}
f(u,s,z)=\frac{u^{z}\zeta(s+z)}{z}
\end{aligned}
\end{equation}
Let us take a rectangular contour $C$ as,
\begin{enumerate}[label=\roman*)]
    \item The vertical line from $2-iT$ to $2+iT$,
    \item The horizontal line from $2+iT$ to $1-\sigma-\eta+iT$,
    \item The vertical line from $1-\sigma-\eta+iT$ to $1-\sigma-\eta-iT$,  
    \item The horizontal line from $1-\sigma-\eta-iT$ to $2-iT$.
\end{enumerate}
where $0<\eta< \min(\sigma,1-\sigma)$ (see Figure \ref{Fig1}). 
\begin{figure}[h]
    \centering
     \includegraphics[width=0.45\textwidth]{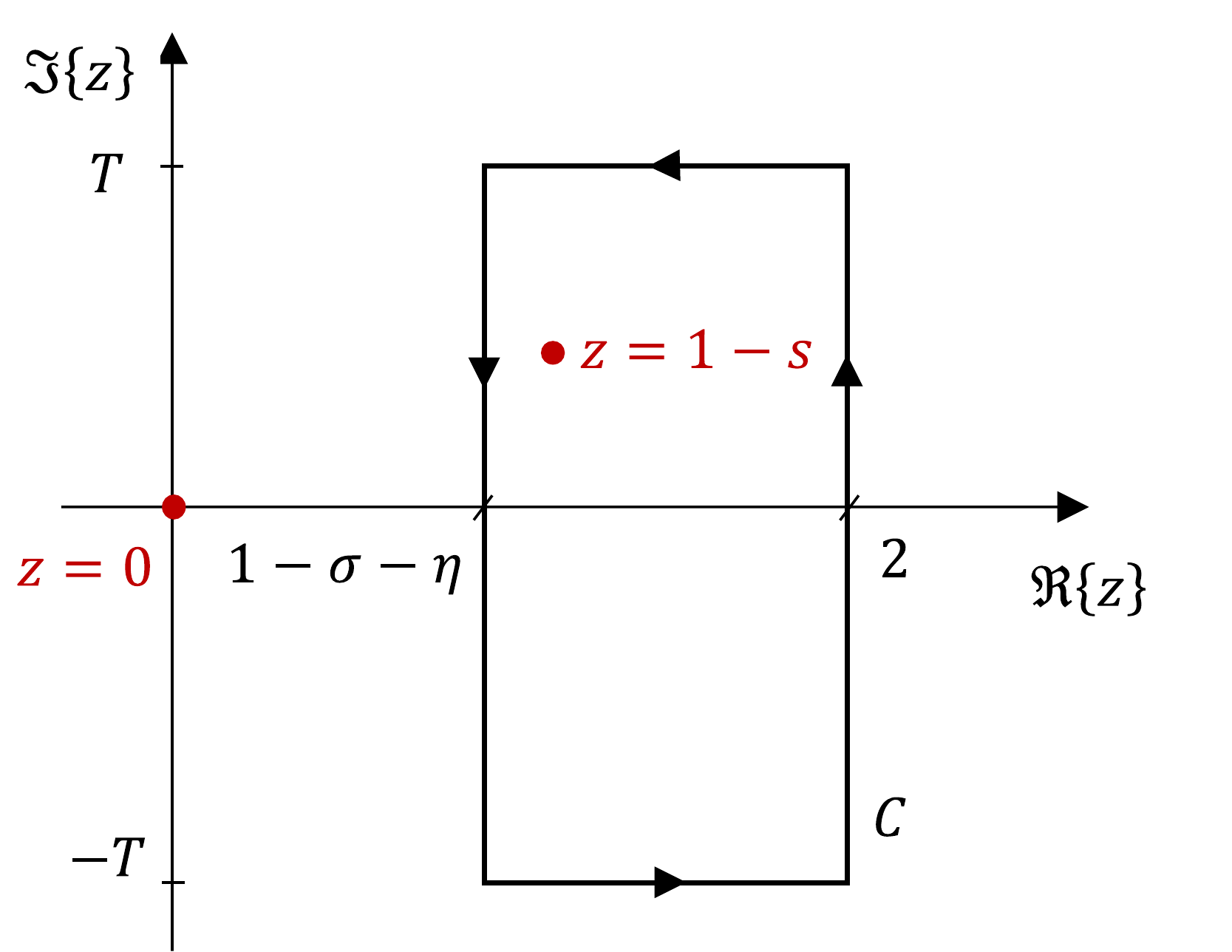}
    \caption{Contour $C$ used to shift the line of integration from $\Re(z)=2$ to $\Re(z)=1-\sigma-\eta$.}
    \label{Fig1}
\end{figure}

Note that only the simple pole at $z=1-s$ of $f(u,s,z)$ lies inside the contour $C$. 

By using analytic continuation of the zeta function and applying Cauchy residue theorem \cite{stein2003}, we get,
\begin{equation}
\begin{aligned}
\label{eq116}
\ointctrclockwise_C f(u,s,z)\,dz&=\int_{2-iT}^{2+iT}f(u,s,z)\,dz+\int_{2+iT}^{1-\sigma-\eta+iT}f(u,s,z)\,dz+\int_{1-\sigma-\eta+iT}^{1-\sigma-\eta-iT}f(u,s,z)\,dz+\int_{1-\sigma-\eta-iT}^{2-iT}f(u,s,z)\,dz\\&=2\pi i \operatorname{Res}_{z=1-s} f(u,s,z)
=2\pi i \left(\frac{u^{1-s}}{1-s}\right)
\end{aligned}
\end{equation}
The second integral can be written as,
\begin{equation}
\begin{aligned}
\label{eq117}
\int_{2+iT}^{1-\sigma-\eta+iT}f(u,s,z)\,dz=\int_{2}^{1-\sigma-\eta}\frac{u^{x+iT}\zeta(s+x+iT)}{x+iT}\,dx
\end{aligned}
\end{equation}
We have,
\begin{equation}
\begin{aligned}
\label{eq118}
|x+iT|\ge T
\end{aligned}
\end{equation}
and from the convexity property of the zeta function \cite{titchmarsh1986}, for large $T$,
\begin{equation}
\begin{aligned}
\label{eq119}
\left|\zeta(s+x+iT)\right|=
\begin{cases} 
O_\varepsilon\left( T^{\eta/2+\varepsilon}\right) &\text{if }1-\sigma-\eta\le x< 1-\sigma \\
O_\varepsilon\left( T^{\varepsilon}\right) &\text{if }1-\sigma\le x\le2  \\
\end{cases}
\end{aligned}
\end{equation}
where $\varepsilon>0$ is an arbitrarily small positive constant, and the implied constant in $O_\varepsilon(\cdot)$ depends on $\varepsilon$.

Therefore,
\begin{equation}
\begin{aligned}
\label{eq120}
\left|\int_{2+iT}^{1-\sigma-\eta+iT}f(u,s,z)\,dz\right|=O_\varepsilon\left( \frac{u^{2}}{T^{1-\eta/2-\varepsilon}}\right)=O_\varepsilon\left( \frac{1}{u^{4-\varepsilon}}\right)
\end{aligned}
\end{equation}
and similarly for the fourth integral.

Substituting in equation (\ref{eq116}), we get
\begin{equation}
\begin{aligned}
\label{eq121}
\frac{1}{2\pi i} \int_{2-iT}^{2+iT}\frac{u^{z}}{z}\zeta(s+z)dz=\frac{1}{2\pi}
\int_{-T}^{T}\frac{u^{1-\sigma-\eta+iy}\zeta(1-\eta+i(t+y))}{1-\sigma-\eta+iy}\,dy+\frac{u^{1-s}}{1-s}+O_\varepsilon\left( \frac{1}{u^{4-\varepsilon}}\right)
\end{aligned}
\end{equation}
Thus, from equations (\ref{eq115}) and (\ref{eq121}), we get
\begin{equation}
\begin{aligned}
\label{eq122}
\sum_{n=1}^{j}\frac{1}{n^s}=\frac{1}{2\pi}\int_{-T}^{T}\frac{u^{1-\sigma-\eta+iy}\zeta(1-\eta+i(t+y))}{1-\sigma-\eta+iy}\,dy+\frac{u^{1-s}}{1-s}+O_\varepsilon\left( \frac{1}{u^{4-\varepsilon}}\right)
\end{aligned}
\end{equation}
From equations (\ref{eq107}) and  (\ref{eq122}), we have
\begin{equation}
\begin{aligned}
\label{eq123}
\zeta(s)=\frac{1}{2\pi}\int_{-u^8}^{u^8}\frac{u^{1-\sigma-\eta+iy}\zeta(1-\eta+i(t+y))}{1-\sigma-\eta+iy}\,dy
+\frac{\bar{B}_1(u)}{u^s}+\frac{s\bar{B}_2(u)}{2u^{s+1}}+\frac{s(s+1)\bar{B}_3(u)}{3!u^{s+2}}+O\left(\frac{\left(1+|t|\right)^4}{u^{\sigma+3}}\right)
\end{aligned}
\end{equation}
Similarly, by replacing $s$ by $(1-\bar{s})$, 
\begin{equation}
\begin{aligned}
\label{eq124}
\zeta(1-\bar{s})=\frac{1}{2\pi}\int_{-u^8}^{u^8}\frac{u^{\sigma-\eta+iy}\zeta(1-\eta+i(t+y))}{\sigma-\eta+iy}\,dy
+\frac{\bar{B}_1(u)}{u^{1-\bar{s}}}+\frac{(1-\bar{s})\bar{B}_2(u)}{2u^{2-\bar{s}}}+\frac{(1-\bar{s})(2-\bar{s})\bar{B}_3(u)}{3!u^{3-\bar{s}}}+O\left(\frac{\left(1+|t|\right)^4}{u^{4-\sigma}}\right)
\end{aligned}
\end{equation}
Assume that $\zeta(s_0)=\zeta(1-\bar{s}_0)=0$, thus, we must have,
\begin{equation}
\begin{aligned}
\label{eq127}
\zeta(s_0)=0=\frac{1}{2\pi}\int_{-u^8}^{u^8}\frac{u^{1-\sigma_0-\eta+iy}\zeta(1-\eta+i(t_0+y))}{1-\sigma_0-\eta+iy}\,dy
+\frac{\bar{B}_1(u)}{u^{s_0}}+\frac{s_0\bar{B}_2(u)}{2u^{s_0+1}}+\frac{s_0(s_0+1)\bar{B}_3(u)}{3!u^{s_0+2}}+O\left(\frac{\left(1+|t_0|\right)^4}{u^{\sigma_0+3}}\right),
\end{aligned}
\end{equation}
\begin{equation}
\begin{aligned}
\label{eq128}
\zeta(1-\bar{s}_0)=0&=\frac{1}{2\pi}\int_{-u^8}^{u^8}\frac{u^{\sigma_0-\eta+iy}\zeta(1-\eta+i(t_0+y))}{\sigma_0-\eta+iy}\,dy
+\frac{\bar{B}_1(u)}{u^{1-\bar{s}_0}}+\frac{(1-\bar{s}_0)\bar{B}_2(u)}{2u^{2-\bar{s}_0}}+\frac{(1-\bar{s}_0)(2-\bar{s}_0)\bar{B}_3(u)}{3!u^{3-\bar{s}_0}}\\&+O\left(\frac{\left(1+|t_0|\right)^4}{u^{4-\sigma_0}}\right)
\end{aligned}
\end{equation}
On the other hand, using equation (\ref{eq4}) with $m=2$, we get
\begin{equation}
\begin{aligned}
\label{eq129}
\zeta(s)&=\sum_{n=1}^j \frac{1}{n^s}-\frac{j^{1-s}}{1-s}-\frac{1}{2 j^s}+\frac{s}{12 j^{s+1}}-\frac{s(s+1)(s+2)}{720 j^{s+3}}+\epsilon(2,j,s),\\
\epsilon(2,j,s)&=-\frac{s(s+1)(s+2)(s+3)(s+4)}{5!}\int_{j}^{\infty} \frac{\bar{B}_5(x)}{x^{s+5}}\,dx
\end{aligned}
\end{equation}
The remainder term can be bounded as,
\begin{equation}
\begin{aligned}
\label{eq130}
|\epsilon(2,j,s)|&\le \frac{\left|s(s+1)(s+2)(s+3)(s+4)\right|\bar{B}_5^{max}}{5!}\int_{j}^{\infty} \frac{1}{x^{\sigma+5}}\,dx&= \frac{\left|s(s+1)(s+2)(s+3)(s+4)\right|\bar{B}_5^{max}}{5!(\sigma+4)j^{\sigma+4}}
\end{aligned}
\end{equation}
where $|\bar{B}_5(x)|<\bar{B}_5^{max}=0.0244582$. 
Therefore,
\begin{equation}
\begin{aligned}
\label{eq131}
\zeta(s)&=\sum_{n=1}^j \frac{1}{n^s}-\frac{j^{1-s}}{1-s}-\frac{1}{2 j^s}+\frac{s}{12 j^{s+1}}-\frac{s(s+1)(s+2)}{720 j^{s+3}}+O\left(\frac{\left(1+|t|\right)^5}{j^{\sigma+4}}\right)
\end{aligned}
\end{equation}
Using Perron's formula as before, for $W=j^8$, we have
\begin{equation}
\begin{aligned}
\label{eq132}
\frac{1}{2\pi i} \int_{c-iW}^{c+iW}\frac{j^z}{z}\zeta(s+z)dz&=\frac{1}{2\pi i} \int_{c-iW}^{c+iW}\frac{j^z}{z}\sum_{n=1}^\infty\frac{1}{n^{s+z}}\,dz \\
&=\frac{1}{2\pi i} \int_{c-iW}^{c+iW}\frac{j^z}{z}\sum_{n=1}^{j-1}\frac{1}{n^{s+z}}\, dz +\frac{1}{2\pi i j^s}\int_{c-iW}^{c+iW}\frac{dz}{z}+\frac{1}{2\pi i} \int_{c-iW}^{c+iW}\frac{j^z}{z}\sum_{n={j+1}}^\infty\frac{1}{n^{s+z}}\, dz \\
\end{aligned}
\end{equation}
Using the following inequality \cite{Apostol1976}, 
\begin{equation}
\label{eq133}
\begin{aligned}
\left|\frac{1}{2\pi i} \int_{c-iW}^{c+iW}\frac{dz}{z}-\frac{1}{2}\right| \le \frac{c}{\pi W}
\end{aligned}
\end{equation}
we get
\begin{equation}
\begin{aligned}
\label{eq134}
\frac{1}{2\pi i} \int_{c-iW}^{c+iW}\frac{j^z}{z}\zeta(s+z)dz&=\frac{1}{2\pi i} \sum_{n=1}^{j-1}\frac{1}{n^s}\int_{c-iW}^{c+iW}\left(\frac{j}{n}\right)^z \frac{dz}{z} +\frac{1}{j^s}\left[\frac{1}{2}+O\left(\frac{c}{\pi W}\right)\right]+\frac{1}{2\pi i} \sum_{n=j+1}^{\infty}\frac{1}{n^s}\int_{c-iW}^{c+iW}\left(\frac{u}{n}\right)^z \frac{dz}{z} 
\end{aligned}
\end{equation}
So, for $n\ge{j}+1$, $n/j\ge (j+1)/j$, then we have
\begin{equation}
\begin{aligned}
\label{eq135}
\left|\frac{1}{2\pi i} \sum_{n=j+1}^{\infty}\frac{1}{n^s}\int_{c-iW}^{c+iW}\left(\frac{j}{n}\right)^z \frac{dz}{z} \right|&\le \sum_{n=j+1}^\infty\frac{1}{n^{\sigma}}\left[\frac{(j/n)^c}{\pi W \log (n/j)}\right]\le \sum_{n=j+1}^\infty\frac{1}{n^{\sigma}}\left[\frac{(j/n)^c}{\pi W \log \left(\frac{j+1}{j}\right)}\right]\\
&\lesssim\frac{j^{c+1}}{\pi W}\sum_{n=j+1}^\infty\frac{1}{n^{\sigma+c}}= O\left(\frac{j^{c+1}}{W}\right)=O\left(\frac{1}{j^{5}}\right)
\end{aligned}
\end{equation}
On the other hand,
\begin{equation}
\begin{aligned}
\label{eq136}
\frac{1}{2\pi i} \sum_{n=1}^{j-1}\frac{1}{n^s}\int_{c-iW}^{c+iW}\left(\frac{j}{n}\right)^z \frac{dz}{z}-
\sum_{n=1}^{j-1}\frac{1}{n^s}&=\sum_{n=1}^{j-1}\frac{1}{n^s}\left[\frac{1}{2\pi i}\int_{c-iW}^{c+iW}\left(\frac{j}{n}\right)^z \frac{dz}{z}-1\right]
\end{aligned}
\end{equation}
For $1\le n\le{j-1}$, $j/n\ge j/(j-1)$, then we have
\begin{equation}
\begin{aligned}
\label{eq137}
\left|\sum_{n=1}^{j-1}\frac{1}{n^s}\left[\frac{1}{2\pi i}\int_{c-iW}^{c+iW}\left(\frac{j}{n}\right)^z \frac{dz}{z}-1\right]\right|
&\le \sum_{n=1}^{j-1}\frac{1}{n^{\sigma}}\left[\frac{(j/n)^c}{\pi W \log \left(\frac{j}{j-1}\right)}\right]\lesssim\frac{(j-1) j^{c}}{\pi W}\sum_{n=1}^{ j-1}\frac{1}{n^{\sigma+c}}\\
&\lesssim\frac{j^{c+1}}{\pi W}\sum_{n=1}^{ j-1}\frac{1}{n^{\sigma+c}}= O\left(\frac{1}{j^{5}}\right)
\end{aligned}
\end{equation}
From equations (\ref{eq134}), (\ref{eq135}) and (\ref{eq137}), we have
\begin{equation}
\begin{aligned}
\label{eq138}
\frac{1}{2\pi i} \int_{2-iW}^{2+iW}\frac{j^z}{z}\zeta(s+z)dz&=\sum_{n=1}^{j-1}\frac{1}{n^s}+\frac{1}{2j^s}+O\left(\frac{1}{j^{5}}\right)&=\sum_{n=1}^{j}\frac{1}{n^s}-\frac{1}{2j^s}+O\left(\frac{1}{j^{5}}\right)
\end{aligned}
\end{equation}
By shifting the contour from $\Re(z)=2$ to $\Re(z)=1-\sigma-\eta$ as before,
\begin{equation}
\begin{aligned}
\label{eq139}
\frac{1}{2\pi i} \int_{2-iW}^{2+iW}\frac{j^{z}}{z}\zeta(s+z)dz=\frac{1}{2\pi}
\int_{-W}^{W}\frac{j^{1-\sigma-\eta+iy}\zeta(1-\eta+i(t+y))}{1-\sigma-\eta+iy}\,dy+\frac{j^{1-s}}{1-s}+O_\varepsilon\left( \frac{1}{j^{4-\varepsilon}}\right)
\end{aligned}
\end{equation}
Thus, from equations (\ref{eq138}) and (\ref{eq139}), we get
\begin{equation}
\begin{aligned}
\label{eq140}
\sum_{n=1}^{j}\frac{1}{n^s}=\frac{1}{2\pi}\int_{-W}^{W}\frac{j^{1-\sigma-\eta+iy}\zeta(1-\eta+i(t+y))}{1-\sigma-\eta+iy}\,dy+\frac{j^{1-s}}{1-s}+\frac{1}{2j^s}+O_\varepsilon\left( \frac{1}{j^{4-\varepsilon}}\right)
\end{aligned}
\end{equation}
From equations (\ref{eq131}) and  (\ref{eq140}), we have
\begin{equation}
\begin{aligned}
\label{eq141}
\zeta(s)=\frac{1}{2\pi}\int_{-j^8}^{j^8}\frac{j^{1-\sigma-\eta+iy}\zeta(1-\eta+i(t+y))}{1-\sigma-\eta+iy}\,dy
+\frac{s}{12 j^{s+1}}-\frac{s(s+1)(s+2)}{720 j^{s+3}}+O\left(\frac{\left(1+|t|\right)^5}{j^{\sigma+4}}\right)
\end{aligned}
\end{equation}
At $s=s_0$,
\begin{equation}
\begin{aligned}
\label{eq142}
\zeta(s_0)=0=\frac{1}{2\pi}\int_{-j^8}^{j^8}\frac{j^{1-\sigma_0-\eta+iy}\zeta(1-\eta+i(t_0+y))}{1-\sigma_0-\eta+iy}\,dy
+\frac{s_0}{12 j^{s_0+1}}-\frac{s_0(s_0+1)(s_0+2)}{720 j^{s_0+3}}+O\left(\frac{\left(1+|t_0|\right)^5}{j^{\sigma_0+4}}\right)
\end{aligned}
\end{equation}
Similarly, by replacing $s_0$ by $(1-\bar{s}_0)$, 
\begin{equation}
\begin{aligned}
\label{eq143}
\zeta(1-\bar{s}_0)=0=\frac{1}{2\pi}\int_{-j^8}^{j^8}\frac{j^{\sigma_0-\eta+iy}\zeta(1-\eta+i(t_0+y))}{\sigma_0-\eta+iy}\,dy
+\frac{1-\bar{s}_0}{12 j^{2-\bar{s}_0}}-\frac{(1-\bar{s}_0)(2-\bar{s}_0)(3-\bar{s}_0)}{720 j^{4-\bar{s}_0}}+O\left(\frac{\left(1+|t_0|\right)^5}{j^{5-\sigma_0}}\right)
\end{aligned}
\end{equation}
Let 
\begin{equation}
\begin{aligned}
\label{eq144}
F(u)&:=\frac{1}{2\pi}\Re\left\{\int_{-u^8}^{u^8}\left[\frac{1}{1-\sigma_0-\eta+iy}-\frac{1}{\sigma_0-\eta+iy}\right]u^{2-\eta+i(t_0+y)}\zeta(1-\eta+i(t_0+y))\,dy\right\},\quad u \ge j_0
\end{aligned}
\end{equation}
From equations (\ref{eq127}) and (\ref{eq128}), for $u\in D$,
\begin{equation}
\begin{aligned}
\label{eq145}
F(u)&=\Re\left\{\frac{[(1-\bar{s}_0)-s_0]\bar{B}_2(u)}{2}+\frac{[(1-\bar{s}_0)(2-\bar{s}_0)-s_0(s_0+1) ]\bar{B}_3(u)}{3!u}+O\left(\frac{\left(1+|t_0|\right)^4}{u^{2}}\right)\right\}\\
&=\frac{(1-2\sigma_0)\bar{B}_2(u)}{2}+\frac{(1-2\sigma_0)\bar{B}_3(u)}{3u}+O\left(\frac{\left(1+|t_0|\right)^4}{u^{2}}\right)
\end{aligned}
\end{equation}
while $F(j)$ can be obtained from equations (\ref{eq142}) and (\ref{eq143}) as,
\begin{equation}
\begin{aligned}
\label{eq146}
F(j)&=\frac{(1-2\sigma_0)}{12}+O\left(\frac{\left(1+|t_0|\right)^2}{j^{2}}\right)
\end{aligned}
\end{equation}
Since the integrand in $F(u)$ is continuous and continuously differentiable with respect to $u$, then using Leibniz rule for differentiation under integral sign \cite{kaplan2003}, we get
\begin{equation}
\begin{aligned}
\label{eq147}
F'(u)&=\frac{1}{2\pi }\Re\left\{\int_{-u^8}^{u^8}\left[\frac{2-\eta+i(t_0+y)}{1-\sigma_0-\eta+iy}-\frac{2-\eta+i(t_0+y)}{\sigma_0-\eta+iy}\right]u^{1-\eta+i(t_0+y)}\zeta(1-\eta+i(t_0+y))\,dy\right.\\
&\qquad\qquad\qquad\left.+8u^7 \left[\frac{1}{1-\sigma_0-\eta+iu^8}-\frac{1}{\sigma_0-\eta+iu^8}\right]u^{2-\eta+i(t_0+u^8)}\zeta(1-\eta+i(t_0+u^8))\right.\\
&\qquad\qquad\qquad\left.+8u^7 \left[\frac{1}{1-\sigma_0-\eta-iu^8}-\frac{1}{\sigma_0-\eta-iu^8}\right]u^{2-\eta+i(t_0-u^8)}\zeta(1-\eta+i(t_0-u^8))\right\}\\
&=\frac{1}{2\pi}\Re\left\{\int_{-u^8}^{u^8}\left[\frac{(s_0+1)}{1-\sigma_0-\eta+iy}-\frac{(2-\bar{s}_0)}{\sigma_0-\eta+iy}\right]u^{1-\eta+i(t_0+y)}\zeta(1-\eta+i(t_0+y))\,dy\right.\\
&\qquad\qquad\qquad\left.+\frac{8u^7(2\sigma_0-1)}{(1-\sigma_0-\eta+iu^8)(\sigma_0-\eta+iu^8)}u^{2-\eta+i(t_0+u^8)}\zeta(1-\eta+i(t_0+u^8))\right.\\
&\qquad\qquad\qquad\left.+ \frac{8u^7(2\sigma_0-1)
}{(1-\sigma_0-\eta-iu^8)(\sigma_0-\eta-iu^8)}u^{2-\eta+i(t_0-u^8)}\zeta(1-\eta+i(t_0-u^8))\right\}\\
\end{aligned}
\end{equation}
From the convexity property of the zeta function,
\begin{equation}
\begin{aligned}
\label{eq148}
|\zeta(1-\eta+i(t_0+u^8))|=O_{\varepsilon} \left(u^{4\eta+\varepsilon}\right)
\end{aligned}
\end{equation}
Thus,
\begin{equation}
\begin{aligned}
\label{eq149}
\left|\frac{8u^7(2\sigma_0-1)}{(1-\sigma_0-\eta+iu^8)(\sigma_0-\eta+iu^8)}u^{2-\eta+i(t_0+u^8)}\zeta(1-\eta+i(t_0+u^8))\right|=O_{\varepsilon} \left(\frac{1}{u^{7-3\eta-\varepsilon}}\right)
\end{aligned}
\end{equation}
and similarly,
\begin{equation}
\begin{aligned}
\label{eq150}
\left|\frac{8u^7(2\sigma_0-1)}{(1-\sigma_0-\eta-iu^8)(\sigma_0-\eta-iu^8)}u^{2-\eta+i(t_0-u^8)}\zeta(1-\eta+i(t_0-u^8))\right|=O_{\varepsilon} \left(\frac{1}{u^{7-3\eta-\varepsilon}}\right)
\end{aligned}
\end{equation}
So,
\begin{equation}
\begin{aligned}
\label{eq151}
F'(u)&=\frac{1}{2\pi}\Re\left\{\int_{-u^8}^{u^8}\left[\frac{(s_0+1)}{1-\sigma_0-\eta+iy}-\frac{(2-\bar{s}_0)}{\sigma_0-\eta+iy}\right]u^{1-\eta+i(t_0+y)}\zeta(1-\eta+i(t_0+y))\,dy \right\}+ O_{\varepsilon} \left(\frac{1}{u^{7-3\eta-\varepsilon}}\right)\\
\end{aligned}
\end{equation}
From equations (\ref{eq127}) and (\ref{eq128}), for $u\in D$
\begin{equation}
\begin{aligned}
\label{eq152}
F'(u)&=O_{\varepsilon} \left(\frac{1}{u^{7-3\eta-\varepsilon}}\right)+\Re\left\{
(s_0+1)\left[u^{s_0}\zeta(s_0)-\bar{B}_1(u)-\frac{s_0\bar{B}_2(u)}{2u}-\frac{s_0(s_0+1)\bar{B}_3(u)}{3!u^{2}}+O\left(\frac{\left(1+|t_0|\right)^4}{u^{3}}\right)\right]\right.\\
&\left.-(2-\bar{s}_0)\left[u^{1-\bar{s}_0}\zeta(1-\bar{s}_0)-\bar{B}_1(u)-\frac{(1-\bar{s}_0)\bar{B}_2(u)}{2u}-\frac{(1-\bar{s}_0)(2-\bar{s}_0)\bar{B}_3(u)}{3!u^{2}}+O\left(\frac{\left(1+|t_0|\right)^4}{u^{3}}\right)\right]\right\}\\
&=(1-2\sigma_0)\bar{B}_1(u)+\frac{(1-2\sigma_0)\bar{B}_2(u)}{u}+O\left(\frac{(1+|t_0|)^2}{u^{2}}\right),\qquad u\in D
\end{aligned}
\end{equation}
On the other hand, $F'(j)$ for $j=\{j_0,j_0+1,\cdots\}$ can be attained from equations (\ref{eq142}) and (\ref{eq143}) as,
\begin{equation}
\begin{aligned}
\label{eq153}
F'(j)&=O_{\varepsilon} \left(\frac{1}{j^{7-3\eta-\varepsilon}}\right)+\Re\left\{
(s_0+1)\left[j^{s_0}\zeta(s_0)-\frac{s_0}{12 j}+\frac{s_0(s_0+1)(s_0+2)}{720 j^{3}}+O\left(\frac{\left(1+|t_0|\right)^5}{j^{4}}\right)\right]\right.\\
&\left.-(2-\bar{s}_0)\left[j^{1-\bar{s}_0}\zeta(1-\bar{s}_0)-\frac{1-\bar{s}_0}{12 j}+\frac{(1-\bar{s}_0)(2-\bar{s}_0)(3-\bar{s}_0)}{720 j^{3}}+O\left(\frac{\left(1+|t_0|\right)^5}{j^{4}}\right)\right]\right\}\\
&=\frac{(1-2\sigma_0)}{6j}+O\left(\frac{(1+|t_0|)^3}{j^{3}}\right)
\end{aligned}
\end{equation}
Since $F'(u)$ is continuous on $[j+a,j+b]$ where $h/j<a<b<1-h/j$ and differentiable on $(j+a,j+b)$, then from the mean value theorem \cite{Trench2003}, we have
\begin{equation}
\begin{aligned}
\label{eq154}
F'\left(j+b\right)-F'\left(j+a\right)=(b-a)F''\left(j+\alpha_j\right),\qquad a<\alpha_j<b
\end{aligned}
\end{equation}
The left-hand side can be obtained from equation (\ref{eq152}) as,
\begin{equation}
\begin{aligned}
\label{eq156}
F'\left(j+b\right)-F'\left(j+a\right)&=(1-2\sigma_0)(b-1/2)+\frac{(1-2\sigma_0)(b^2-b+1/6)}{j+b}+O\left(\frac{(1+|t_0|)^2}{(j+b)^{2}}\right)\\
&-(1-2\sigma_0)(a-1/2)-\frac{(1-2\sigma_0)(a^2-a+1/6)}{j+a}-O\left(\frac{(1+|t_0|)^2}{(j+a)^{2}}\right)\\
&=(1-2\sigma_0)(b-a)+\frac{(1-2\sigma_0)(b-a)(a+b-1)}{j}+O\left(\frac{(1+|t_0|)^2}{j^{2}}\right)
\end{aligned}
\end{equation}
where equation (\ref{eq103}) was used.

Thus, from equation (\ref{eq154}),
\begin{equation}
\begin{aligned}
\label{eq157}
F''\left(j+\alpha_j\right)=(1-2\sigma_0)+\frac{(1-2\sigma_0)(a+b-1)}{j}+O\left(\frac{(1+|t_0|)^2}{j^{2}}\right)
\end{aligned}
\end{equation}
Since this equation holds for all possible values of $a,b$ and $h$ such that $h/j<a<b<1-h/j$, then as $j\to\infty$, the whole interval $(0,1)$ is spanned uniformly and we have,
\begin{equation}
\begin{aligned}
\label{eq158}
\lim_{u\to\infty}F''\left(u\right)=(1-2\sigma_0),\qquad {u}=u-\lfloor u \rfloor \in (0,1)
\end{aligned}
\end{equation}
Now, let us apply the mean value theorem between the points $j$ and $(j+e)$, $h/j<e<1-h/j$,
\begin{equation}
\begin{aligned}
\label{eq159}
F'\left(j+e\right)-F'\left(j\right)=eF''\left(j+\beta_j\right),\qquad 0<\beta_j<e
\end{aligned}
\end{equation}
Taking the limit as $j\to\infty$,
\begin{equation}
\begin{aligned}
\label{eq160}
\lim_{j\to\infty}\left[F'\left(j+e\right)-F'\left(j\right)\right]=e \lim_{j\to\infty} F''\left(j+\beta_j\right)
\end{aligned}
\end{equation}
From equations (\ref{eq152}),(\ref{eq153}) and (\ref{eq158}),
\begin{equation}
\begin{aligned}
\label{eq161}
&\lim_{j\to\infty}\left[(1-2\sigma_0)\left(e-1/2\right)+\frac{(1-2\sigma_0)\bar{B}_2(e)}{j+e}+O\left(\frac{(1+|t_0|)^2}{(j+e)^{2}}\right)-\frac{(1-2\sigma_0)}{6j}-O\left(\frac{(1+|t_0|)^3}{j^{3}}\right)\right]=e(1-2\sigma_0)
\end{aligned}
\end{equation}
That is, 
\begin{equation}
\begin{aligned}
\label{eq162}
-\frac{1}{2}(1-2\sigma_0)=0 \Rightarrow \sigma_0=1/2
\end{aligned}
\end{equation}
Therefore, all the non-trivial zeros of the zeta function must lie on the critical line, $\Re(s)=1/2$, concluding the proof of the Riemann hypothesis.
\end{proof}

\bibliographystyle{unsrt}
\bibliography{Perron_int_partial56111723}
\end{document}